\documentclass[11pt]{amsart}
\begin{document}

\title{Addendum to ``Structure of seeds in generalized cluster algebras''}
\address{\noindent Graduate School of Mathematics, Nagoya University, 
Chikusa-ku, Nagoya,
464-8604, Japan}
\email{nakanisi@math.nagoya-u.ac.jp}
\author{Tomoki Nakanishi}
\date{}
\maketitle

\begin{abstract}
This is Addendum to 
``Structure of seeds in generalized cluster algebras'', Pacific J. Math. {277} (2015), 201--218.
We extend the  class of generalized cluster algebras studied therein
to embrace examples in some  applications.
\end{abstract}

The definition of the class of generalized cluster algebras
in \cite{Nakanishi14a} turned out to be too restrictive
for some  applications (e.g., \cite{Gyoda22}).
However, it is easily extended as follows.
Let $ \mathbb{P}$ be the semifield in Section 2A.
In (2.1)
we replace each data $z_{i,s}\in \mathbb{P}$ 
with $z_{i,s}\in \mathbb{Z}_{\geq 0}\mathbb{P}$.
For $z_{i,s}\neq 0$,
let $z_{i,s}=\sum_{j} m_j p_j$ ($m_j\in \mathbb{Z}_{> 0}$, $p_j\in \mathbb{P}$),
where  the multiplication by $m_j $ is in $\mathbb{Z}_{\geq 0}\mathbb{P}$.
Then, 
we set $\tilde z_{i,s}=\bigoplus_{j} m_j p_j\in \mathbb{P}$, where the multiplication by $m_j $ is in the semifield
$\mathbb{P}$.
For  $z_{i,s}= 0$, we formally set $\tilde z_{i,s}=0$, though this is not an element in $\mathbb{P}$.
Then, we replace the formulas (2.5) and (2.6) with
\begin{align*}
{y}'_i&=
\begin{cases}
\displaystyle
{y}_k^{-1}
 & i=k\\
 \displaystyle
{y}_i
\Biggl(
{y}_k^{[\varepsilon {b}_{ki}]_+} 
\Biggr)^{d_k}
\Biggl(
\bigoplus_{s=0}^{d_k}
\tilde z_{k,s} {y}_k^{\varepsilon s}
\Biggr)^{-{b}_{ki}}
& i\neq k,\\
\end{cases}
\\
x'_i&=
\begin{cases}
\displaystyle
x_k^{-1}
\Biggl(\prod_{j=1}^n 
x_j^{[-\varepsilon {b}_{jk}]_+}
\Biggr)^{d_k}
\frac{\displaystyle
\sum_{s=0}^{d_k}
 z_{k,s} \hat{y}_k^{\varepsilon s}
}
{\displaystyle
\bigoplus_{s=0}^{d_k}
\tilde  z_{k,s} {y}_k^{\varepsilon s}
}
 & i=k\\
x_i& i\neq k,\\
\end{cases}
\end{align*}
where, if $\tilde z_{k,s}=0$ appears  in the sum in $\mathbb{P}$,  we simply ignore it.

If we set $\mathbb{P}$  to be the trivial semifield,
it embraces the examples in \cite{Gyoda22}.
If there is some $z_{i,s}\not\in \mathbb{P}$,
it is no longer a generalized cluster algebra  in the sense of \cite{Chekhov11};
accordingly, the correspondence given in Section 2A after Proposition 2.5 does not hold.
However, all  the other propositions and theorems  remain valid with the obvious modification,
namely, replacing $F_i^t\vert_{\mathbb{P}}(\mathbf{y},\mathbf{z})$ with
$F_i^t\vert_{\mathbb{P}}(\mathbf{y},\tilde{\mathbf{z}})$ in Theorems 3.22 and 3.23.
We thank Yasuaki Gyoda and Kodai Matsushita for discussion.

\bibliographystyle{amsalpha}
\bibliography{../../biblist/biblist.bib}

\providecommand{\bysame}{\leavevmode\hbox to3em{\hrulefill}\thinspace}
\providecommand{\MR}{\relax\ifhmode\unskip\space\fi MR }
\providecommand{\MRhref}[2]{%
  \href{http://www.ams.org/mathscinet-getitem?mr=#1}{#2}
}
\providecommand{\href}[2]{#2}
\begin{thebibliography}{Nak15}

\bibitem[CS14]{Chekhov11}
L.~Chekhov and M.~Shapiro, \emph{Teichm\"uller spaces of {R}iemann surfaces
  with orbifold points of arbitrary order and cluster variables}, Int. Math.
  Res. Notices \textbf{2014} (2014), 2746--2772; arXiv:1111.3963 [math--ph].

\bibitem[GM23]{Gyoda22}
Y.~Gyoda and K.~Matsushita, \emph{Generalization of {M}arkov {D}iophantine
  equation via generalized cluster algebra}, Electron. J. Combin. \textbf{30}
  (2023), P4.10; arXiv: arXiv:2201.10919 [math.NT].

\bibitem[Nak15]{Nakanishi14a}
T.~Nakanishi, \emph{Structure of seeds in generalized cluster algebras},
  Pacific J. Math. \textbf{277} (2015), 201--218; arXiv:1409.5967 [math.RA].

\end{thebibliography}

\end{document}